\RequirePackage[leqno]{amsmath}
\documentclass[leqno]{amsart}
\usepackage{color}
\usepackage{amsfonts}
\usepackage{enumitem}
\setlist[itemize]{topsep=0pt,after=\vspace{1.5\baselineskip}}
 \usepackage[colorlinks=true]{hyperref}
\usepackage{empheq}
\usepackage[T1]{fontenc}
\usepackage[latin1]{inputenc}

\def\R{\mathbb R}  
\newtheorem{theorem}{Theorem}[section]

\newtheorem{lemma}[theorem]{Lemma}

\newtheorem{remark}{Remark}

\title[Boundedness in a general chemotaxis system with logistic source] 
      {Boundedness in a fully parabolic chemotaxis system with nonlinear diffusion and sensitivity, and logistic source}
\author[G. Viglialoro and M. Marras]{}
\subjclass[2010]{35K55, 35Q92, 35A01, 92C17.}
\keywords{Nonlinear parabolic systems, chemotaxis, logistic source, global existence, boundedness.\\
\textit{$^*$Corresponding author}: giuseppe.viglialoro@unica.it}
\begin{document}
\maketitle
\centerline{\scshape  M. Marras$^1$ \and G. Viglialoro$^{1,*}$}
\medskip
{\footnotesize
 \centerline{$^1$Dipartimento di Matematica e Informatica}
 \centerline{Universit\`{a} di Cagliari}
 \centerline{V. le Merello 92, 09123. Cagliari (Italy)}
 \medskip
 \medskip
}

\bigskip
\begin{abstract}
In this paper we study the zero-flux chemotaxis-system 
\begin{equation*}
\begin{cases}
u_{ t}=\nabla \cdot ((u+1)^{m-1} \nabla u-(u+1)^\alpha \chi(v)\nabla v)+ku-\mu u^2 & x\in \Omega, t>0, \\
v_{t}=\Delta v-vu & x\in \Omega,  t>0,\\
\end{cases}
\end{equation*}
$\Omega$ being a bounded and smooth domain of $\mathbb{R}^n$, $n\geq 1$, and where $m,k \in \R$, $\mu>0$ and $\alpha < \frac{m+1}{2}$.  For any $v\geq 0$ the  chemotactic sensitivity function is assumed to behave as the prototype $\chi(v)=\frac{\chi_0}{(1+av)^2}$, with $a\geq 0$ and $\chi_0>0$. We prove that for nonnegative and sufficiently regular initial data $u(x,0)$ and $v(x,0),$ the corresponding initial-boundary value problem admits a global bounded classical solution provided $\mu$ is large enough.
\end{abstract}
\section{Introduction and motivations}\label{IntroductionSection} 
The expression \textit{chemotaxis} indicates the movement of cells occupying a space, which are stimulated by a chemical signal produced by a substance therein inhomogeneously distributed. Different studies and experiments on bacteria show how they direct their natural motion and change randomly their course precisely depending on the intensity of the chemical gradient stimulus.

In 1971, Keller and Segel (see \cite{Keller-1971-TBC}) proposed a model for the description of  the traveling band behavior of bacteria due to the chemotactic response; its mathematical formulation defined in an interval $I\subset \R$ reads:
\begin{equation}\label{Keller-Segel-1971_Formulation}
\left\{ \begin{array}{ll}
u_{t}=  u_{xx}-\chi (uv^{-1}v_x)_{x}& x\in I, t>0, \\
v_{t}= \varepsilon v_{xx}-uf(v) &x\in I, t>0.\end{array} \right. 
\end{equation}
In this model, the distribution of the cells  and the concentration of the chemical signal in a point $x$ of $I$ and at an instant $t$ of time are, respectively, identified with the functions $u=u(x,t)$ and $v=v(x,t)$. In addition, $\varepsilon>0$ represents the diffusion coefficient of the chemical substance, $\chi$ is called chemosensitivity and essentially measures the drift velocity of the bacteria through the direction of the gradient of concentration of the chemoattractant and $f (v)$ denotes a kinetic function describing the chemical reaction between bacteria and the chemical. In terms of the value of the constant $\varepsilon$ and the expression of  the function $f$, some results concerning the existence of traveling wave solutions to system \eqref{Keller-Segel-1971_Formulation} have been established, under proper boundary and initial conditions. For instance, the same \cite{Keller-1971-TBC} corresponds to the limit case  $\varepsilon=0$ and $f (v) = \alpha> 0$, \cite{NagaiIkeda_Travelin} to $\varepsilon>0$ and $f (v) = \alpha> 0$ and \cite{LiWangTravelling} to  $\varepsilon>0$ and $f (v) = \alpha v> 0$, with $\alpha>0.$

A possible extension to higher dimensions of  \eqref{Keller-Segel-1971_Formulation} is described by this initial-boundary value problem 
\begin{equation}\label{Keller-Segel-CosnumptionTao}
\left\{ \begin{array}{ll}
u_{t}=  \Delta u - \chi \nabla ( u \cdot \nabla v)& x\in \Omega, t>0, \\
v_{t}=\Delta v-uv & x\in \Omega,  t>0,\\
\frac{\partial u}{\partial \nu}=\frac{\partial v}{\partial \nu}=0 & x\in \partial \Omega, t>0, \\
u(x,0)=u_{0}(x) \geq 0 \,\,\textrm{and}\,\,\,v(x,0)=v_{0}(x) \geq 0& x\in  \Omega,
\end{array} \right. 
\end{equation}
where $\chi>0$, $\Omega \subset \R^n$, with $n\geq 2$, is a bounded domain with smooth boundary,  and $u_0(x)=u(x,0)$ and  $v_0(x)=v(x,0)$ are the initial cells distribution and chemical concentration. Moreover, since $\frac{\partial}{\partial \nu}$ indicates the outward normal derivative on $\partial \Omega$,  with zero-flux boundary conditions on both $u$ and $v$, \eqref{Keller-Segel-CosnumptionTao} describes the dynamic of a cells population in response to a chemical substance which mutually interact in a totally insulated domain.

Let us observe that for positive cells and chemical distributions, the term $-uv$ in the second equation of \eqref{Keller-Segel-CosnumptionTao} shows that in such a model the signal is progressively consumed by cells; this, naturally, suggests that $v$ remains bounded through the time. This situation is far to be similar in the classical Keller-Segel model (see \cite{K-S-1970}), where $-uv$ reads $-v+u$ and, hence, the corresponding formulation 
\begin{equation}\label{Keller-Segel-calssic}
\left\{ \begin{array}{ll}
u_{t}=  \Delta u - \chi \nabla ( u \cdot \nabla v)& x\in \Omega, t>0, \\
v_{t}=\Delta v-v+u & x\in \Omega,  t>0,\\
\frac{\partial u}{\partial \nu}=\frac{\partial v}{\partial \nu}=0 & x\in \partial \Omega, t>0, \\
u(x,0)=u_{0}(x) \geq 0 \,\,\textrm{and}\,\,\,v(x,0)=v_{0}(x) \geq 0& x\in  \Omega,
\end{array} \right. 
\end{equation}
manifests how an increasing of the cells favors a production of the signal, so that no bound for $v$ is a priori expected. Thereafter,  even though they are deeply connected, the two models \eqref{Keller-Segel-CosnumptionTao} and \eqref{Keller-Segel-calssic} present different properties. In particular, for system \eqref{Keller-Segel-calssic}, which has been widely discussed from many authors, it is known that in a one-dimensional domain its solutions are global and uniformly bounded in time (see \cite{OsYagUnidim}), while in the $n$-dimensional context, with $n\geq 2$, unbounded solutions idealizing the so called \textit{chemotactic collapse} (an uncontrolled gathering of cells in proximity of some zones of the space occurring in a certain blow-up time) have been detected e.g. in \cite{HorstWang} and \cite{W}. In accordance with this, estimates from below for the blow-up time of unbounded solutions to \eqref{Keller-Segel-calssic} are explicitly derived in \cite{PSong}. 

Now, in order to better contextualize this present investigation in the frame of the existent literature, we precise that an analysis of the previous last contributions highlights as,  in higher dimensions ($n\geq 2$) and under suitable smallness assumptions on the initial data $u_0$ and $v_0$, the solution to model \eqref{Keller-Segel-calssic} is global and bounded while there exist blow-up solutions to model \eqref{Keller-Segel-calssic} for large initial data $u_0$. On the opposite side,  it has been shown that the global existence or blow-up
of solutions to model \eqref{Keller-Segel-CosnumptionTao} is independent of the initial data $u_0$. Specifically, Tao (see \cite{TaoBoun}) proves that, for sufficiently regular $u_0$ and $v_0$, if 
\begin{equation}\label{Conditionv0chiTao}  
0 < \chi \lVert v_0 \rVert_{L^\infty(\Omega)}\leq \frac{1}{6(n+1)},
\end{equation}
then the corresponding initial-boundary value problem \eqref{Keller-Segel-CosnumptionTao} possesses a unique global solution that is uniformly bounded, remaining indeed still open the question  whether there exist  blow-up solutions to the same system for large initial data $v_0$ or large chemotactic parameter $\chi$ not complying with \eqref{Conditionv0chiTao}. Moreover, continuing on the state of the art of the chemotaxis-consumption model  \eqref{Keller-Segel-CosnumptionTao}, for the three-dimensional setting, in \cite{TaoWinkConsumptionEventual} weak solutions that become smooth after some time are constructed. Further, by interpreting the same system  as the special case of the general coupled chemotaxis-fluid model, proposed by Goldstein in \cite{TuvalGoldsteinEtAl}, where exactly the fluid does not give any direct contribution, in \cite{WinklerN-Sto_CPDE} the existence of global classical and weak solutions for $n=2$ and $n=3$, respectively, is discussed while \cite{WinklerN-Sto_2d} deals with the stabilization properties of these two-dimensional solutions (we also cite \cite{WinklerTensorialConsumption} and \cite{ZhangTensorialConsumption} for existence results to systems close to the same \eqref{Keller-Segel-CosnumptionTao} but involving matrix-valued sensitivities). 

Exactly in order to contrast undesired blow-up singularities, which as mentioned above may emerge in both models \eqref{Keller-Segel-CosnumptionTao} and \eqref{Keller-Segel-calssic}, more complete formulations to these systems with nontrivial sources have been considered; precisely, a complementation of these models through largely used \textit{logistic-type} effects seems totally natural (see also \cite{Marras-Viglia-ComptesGradient} and \cite{VigliaGradTermDiffIntEqua} for another expression for the source). For instance, for the system  
\begin{equation}\label{sys Johannes}
\begin{cases}
 u_{t}=  \Delta u -  \chi \nabla \cdot (u \nabla v) + g(u) & x\,\in \, \Omega, \,\,\,t>0, \\
 \tau v_{t}=  \Delta v- v+ u  & x\,\in \, \Omega, \,\,\,t>0,          
 \end{cases}
\end{equation}
defined in a convex smooth and bounded domain $\Omega$ of $\mathbb{R}^n$, $n\geq 1$ and endowed with homogeneous Neumann boundary conditions,  in \cite{Lankeit} the existence of global weak solutions for $\chi=\tau=1$ and $g(u)=ku-\mu u^2$,  for $k \in \R$ and $\mu$ positive constant (the classical logistic source), is established for any nonnegative and sufficient regular  initial data $(u_0,v_0)$ and arbitrarily small values of $\mu>0$;  moreover, if $n=3$, these solutions become classical after some time and provided that $k$ is not too large. On the other hand, under the same assumptions on the domain and the boundary conditions, for a source $g$ generalizing the logistic source and verifying $g(0) \geq 0$ and $g(s)\leq k -\mu s^2$, for $s \geq 0$, $k \geq 0$ and $\mu, \tau$ positive and $\chi \in \R$,  in \cite{W0} the author proves that if  $\mu$ is big enough, for all sufficiently smooth and nonnegative initial data $u_0$ and $v_0$, system \eqref{sys Johannes} possesses a unique bounded and global-in-time classical solution.  Additionally, for the same problem \eqref{sys Johannes}, also defined in a convex smooth and bounded domain $\Omega$ of $\mathbb{R}^n$, $n\geq 1$, but with source term $g$ such that $ -c_0(s+s^\alpha)\leq g(s)\leq a -b s^\alpha$, for $s \geq 0$, and  with some $\alpha>1$, $a\ge 0$ and $\chi, b ,c_0>0$, global existence of very weak solutions, as well their boundedness properties and long time behavior are discussed in \cite{ViglialoroVeryWeak}, \cite{ViglialoroBoundnessVeryWeak} and also \cite{ViglialoroWolleyDCDS}.  Finally, for the sake of completeness, it is worth to precise that even though in the logistic source the term $-\mu u^2$, with  $\mu>0$,  corresponds to a  death rate of the cells distribution and generally contrasts blow-up phenomena,  in \cite{WinDespiteLogistic} is shown that under radially symmetric assumptions there exist initial data such that the corresponding solution of systems type \eqref{sys Johannes} blows up  (we also refer to \cite{MarrasVernierVigliaWithm} for techniques used to estimate the blow up time of unbounded solutions to system related to \eqref{sys Johannes}).

To the best of our knowledge, deserving results in the direction of this present investigation which are tied to \eqref{Keller-Segel-CosnumptionTao} under a perturbation for the evolution of $u$ through a logistic source, are the following. In a bounded and smooth domain of $\R^n$, $n\geq 1$, and under homogeneous Neumann boundary conditions, the system  
\begin{equation}\label{sysGeneralForConsumptionAndLogisticJohannes}
\begin{cases}
 u_{t}=  \nabla \cdot (D(u)\nabla u) - \nabla \cdot (u \chi(v)\nabla v) + g(u) & x\,\in \, \Omega, \,\,\,t>0, \\
  v_{t}=  \Delta v- uf(v)  & x\,\in \, \Omega, \,\,\,t>0,          
 \end{cases}
\end{equation}
\begin{enumerate} [label=\roman*)]
\item  \label{BaghaeiaKhelghatibItems} for $D(u)= \delta >0$, $\chi(v)=\chi>0$, $g(u)\leq k u-\mu u^\gamma$, with $g(0)\geq 0$, $k\geq 0$, $\mu>0$ and $\gamma>1$, and sufficiently regular initial data $(u_0,v_0),$ admits for suitable small $\chi \lVert v_0\rVert_{L^\infty(\Omega)}$ a global and bounded classical solution (see \cite{BaghaeiaKhelghatib}); 
\item  \label{ZhengMuItems} for $D(u)\equiv 1$, $\chi(v)\leq \frac{\chi_0}{(1+a s)^b}$, with $\chi_0, a,b>0$, $g(u)\leq  k u-\mu u^2$, with $g(0)=0$ and $k,\mu>0$,  and sufficiently regular initial data $(u_0,v_0,)$ admits for $\chi_0 \lVert v_0\rVert_{L^\infty(\Omega)}\leq \frac{1}{6(n+1)}$ a global and bounded classical solution (see \cite{ZhengMuConsumption});
\item  \label{LankeitConsumptionItem} for $D(u)\equiv 1$, $\chi(v)=\chi>0$, $k \in \R$, $g(u)= ku-\mu u^2$, and sufficiently regular initial data $(u_0,v_0),$ admits for $\mu$ larger compared to $\chi \lVert v_0\rVert_{L^\infty(\Omega)}$ a global and bounded classical solution and a weak one for arbitrary $\mu>0$ (see \cite{LankeitWangConsumptLogistic}). 
\end{enumerate}
\begin{remark}\label{RemarckChiV_0}
As shown in the corresponding proofs, we conclude this section observing that  the contribution summarized in \ref{LankeitConsumptionItem} explicitly provides conditions which connects $v_0$ to the coefficient $\mu$ of the logistic term $g$; this does not hold for the items \ref{BaghaeiaKhelghatibItems} and \ref{ZhengMuItems}.  Nevertheless,  the term $\chi \lVert v_0\rVert_{L^\infty(\Omega)}$ appears in these same works (as well in the mentioned assumption \eqref{Conditionv0chiTao}  of \cite{TaoBoun}), so that it seems a very proper quantity which coherently characterizes the nature of models type  \eqref{sysGeneralForConsumptionAndLogisticJohannes}; consequently it will play a crucial role also in our main result.
\end{remark}
\section{Main result and structure of the paper}\label{StructureMainReusltSection}
In agreement with all of the above, this paper is dedicated to the following problem
\begin{equation}\label{problem}
\begin{cases}
u_{ t}=\nabla \cdot ((u+1)^{m-1} \nabla u-(u+1)^\alpha \chi(v)\nabla v)+ku-\mu u^2 & x\in \Omega, t>0, \\
v_{t}=\Delta v-vu & x\in \Omega,  t>0,\\
\frac{\partial u}{\partial \nu}=\frac{\partial v}{\partial \nu}=0 & x\in \partial \Omega, t>0, \\
u(x,0)=u_{0}(x) \geq 0 \,\,\textrm{and}\,\,\,v(x,0)=v_{0}(x) \geq 0& x\in  \Omega,
\end{cases}
\end{equation}
defined in a bounded and smooth domain $\Omega$ of $\mathbb{R}^n$, $n\geq 1$, where $(u_{0}, v_{0})$ is a pair of nonnegative functions from $(W^{1,r}(\Omega))^2,$ for some $r>\max\{n,2\}$, and $m,k,\mu \in \R$ with $\mu>0$. Moreover, we assume that 
\begin{equation}\label{conditiononAlpha}
\alpha < \frac{m+1}{2},
\end{equation}
and that the function $\chi$ generalizes the standard \textit{chemotactic sensitivity} 
\begin{equation}\label{StadardChemotactic}  
\chi(v)= \frac{\chi_0}{(1+av)^2} \quad \textrm{with some }\quad \chi_0>0 \quad \textrm{and}\quad a\geq 0;
\end{equation} 
exactly,  $\chi\in C^2([0,\infty))$ and satisfies for any $v\geq 0$ this growth condition
\begin{equation}\label{GrowthConditionOnS}  
\chi(v)\leq \frac{\chi_0}{(1+av)^b} \quad \textrm{with some }\quad \chi_0>0 \quad \textrm{and}\quad a\geq 0 \quad \textrm{and}\quad b>0.
\end{equation} 
Specifically, the aim of the present article is to prove the existence of global and bounded classical solutions to \eqref{problem} under some largeness assumption on $\mu$ (whose value will depend on $n, m$ and $\alpha$) with respect to some combination of powers of $\chi_0\lVert v_0\rVert_{L^\infty(\Omega)}$, precisely in accordance to Remark \ref{RemarckChiV_0} of the previous section.
\begin{remark}
The problem studied in this paper generalizes some of the aforementioned examples. For instance, starting from the case in which no logistic source affects the system (i.e. $k=\mu=0$), in \eqref{problem} the limit values $m=\alpha=1$ and $\chi(v)=\chi >0$ recover \eqref{Keller-Segel-CosnumptionTao}; further,  in presence of logistic perturbations, the models discussed in items  \ref{BaghaeiaKhelghatibItems},  \ref{ZhengMuItems} and \ref{LankeitConsumptionItem} are also easily deducible from \eqref{problem} through evident choices of its data. 

Additionally, we also observe that the first equation of \eqref{problem} can be equivalent rewritten as 
\begin{equation*}
u_{ t}=\nabla \cdot ((u+1)^{m-1} \nabla u-(u+1)^\alpha\nabla \Theta (v) )+ku-\mu u^2 \quad  x\in \Omega, t>0, 
\end{equation*}
where in view of assumption \eqref{GrowthConditionOnS} we have that $\Theta (v)=-\int_v^\infty \chi(s)ds$ is finite. In this sense, it could be checked that the main result herein applies whenever $\Theta$ is smooth and nondecreasing on $[0,\infty)$ and satisfies
\begin{equation*}
\Theta (s)=\frac{s^\beta}{1+as^\beta}\quad a\geq 0,\,\beta>0 \quad \textrm{and}\quad s \geq 0.
\end{equation*}
This last expression for $\Theta$ includes also the cases below, which are biologically coherent (see \cite{Murray1} and \cite{Murray2}),
\begin{equation*}
\Theta (s)=\frac{s}{1+as}\quad \textrm{and}\quad\Theta (s)=\frac{s^2}{1+as^2} \quad \textrm{and}\quad s \geq 0,
\end{equation*}
the first one behaving actually as our prototype \eqref{StadardChemotactic}.
\end{remark}
The rest of the paper is organized as follows. Once  in $\S$\ref{PreliminariesSection} some preparatory and well-known preliminaries are given, $\S$\ref{ExistenceSolutionRegularizingSection} is focused on the derivation of a result concerning local-in-time existence of a classical solution $(u,v)$ to system \eqref{problem} and on some crucial properties tied to the $u$- and $v$-components.  Thereafter, in $\S$\ref{EstimatesAndProofSection}, we define a sort of energy for such a local solution, which for some suitable $p>1$ is defined as $\Phi (t):=\int_\Omega (u+1)^p+\chi_0^{2p}\int_\Omega |\nabla v|^{2p}$. In this way, by relying on general functional and algebraic relations,  we establish that this energy satisfies a differential inequality, under a certain initial condition, whose right hand side is a power function which possesses a positive root. Subsequently, by means of a comparison principle, we provide a local-in-time independent bound for $\Phi (t)$, and hence in particular for $u$ in $L^p(\Omega)$ and $\nabla v$ in $L^{2p}(\Omega)$. Finally, also in the same $\S$\ref{EstimatesAndProofSection}, an application of a general and extensively used boundedness result for parabolic equations allows us to deduce from  these gained estimates that the local classical solution is actually bounded and, then, also global. This represents exactly our main assertion: 
\begin{theorem}\label{MainTheorem}
Let $\Omega$ be a smooth and bounded domain of $\mathbb{R}^n$, with $n\geq 1$. For given $m,k\in \R$, and $\mu$ positive, let us assume that  $\chi\in C^2([0,\infty))$ satisfies relation \eqref{GrowthConditionOnS} for some $\alpha$ as in  \eqref{conditiononAlpha}. Then for any couple of nonnegative functions $(u_0,v_0)\in(W^{1,r}(\Omega))^2$, with $r>\max\{n,2\}$, it is possible to find two positive constants $K_1 = K_1(n,m,\alpha)$ and $K_2 = K_2(n,m,\alpha)$ such that if 
\begin{equation}\label{LargenessAssumptionMu}
\mu > K_1(n,m,\alpha)\lVert \chi_0 v_0\rVert^\frac{2}{n}_{L^\infty(\Omega)}+K_2(n,m,\alpha)\lVert \chi_0 v_0 \rVert^{2n}_{L^\infty(\Omega)},
\end{equation}
problem \eqref{problem} admits a unique global classical solution $(u, v)$ which is uniformly bounded; more precisely, there exists a positive constant $C$ such that
\begin{equation}\label{GlobalBoundednessuinftyvW1infty}
\lVert u(\cdot,t) \rVert_{L^\infty(\Omega)}+\lVert v(\cdot,t) \rVert_{W^{1,r}(\Omega)}\leq C \quad \textrm{for all} \quad t\in(0,\infty).
\end{equation}
\end{theorem}
%
\section{Some preparatory tools}\label{PreliminariesSection}
The following Lemmas are used through the paper to prove the main theorem. In particular, we mainly summarize some general functional inequalities, other proper technical results, and we close the section by adjusting some parameters which are necessary in our logical steps. 
\begin{lemma}\label{InequalityTipoPoincarLemma} 
Let $\Omega$ be a bounded and smooth domain of $\R^n$, $n\geq 1,$ and $q\geq 1$. For all $f\in C^{2}(\Omega)$, we have
\begin{equation}\label{InequalityLaplacian}
(\Delta f)^2 \leq n |D^2 f|^2,
\end{equation}
\begin{equation}\label{InequalityGradienHessian}
\lvert D^2f \nabla f \rvert^2\leq |D^2 f|^2\lvert \nabla f \rvert^2,
\end{equation}
and for all $f\in C^{2}(\bar{\Omega})$ satisfying $f\frac{\partial f}{\partial \nu}=0$ on $\partial \Omega$,
\begin{equation}\label{InequalityHessian}
\lVert  \nabla f\rVert^{2q+2}_{L^{2q+2}(\Omega)} \leq 2 (4q^2+n)\lVert f \rVert^2_{L^\infty(\bar{\Omega})} \lVert \vert \nabla f\vert^{q-1}  D^2 f\rVert^2_{L^2(\Omega)},
\end{equation}
where $D^2f$ represents the Hessian matrix of $f$ and $\vert D^2f \rvert^2=\sum\limits_{i,j=1}^{n}f_{x_ix_j}^2.$
\begin{proof}
Straightforward calculations infer
\[
\begin{split}
(\Delta f)^2=\Big(\sum_{i=1}^{n}f_{x_i x_i}\Big)^2&=\sum_{i,j=1}^{n}f_{x_i x_i}f_{x_j x_j}   \leq \sum_{i,j=1}^{n}\Big(\frac{1}{2}f^2_{x_i x_i}+ \frac{1}{2}f^2_{x_j x_j}\Big)\\ & = n \sum_{i=1}^{n}f^2_{x_i x_i}\leq n  \sum_{i,j=1}^{n}f^2_{x_i x_j}=n\lvert D^2 f\rvert^2,
\end{split}
 \]
 and that 
 \[
 \begin{split}
\lvert D^2 f \nabla f\rvert^2 &=\sum_{i,j=1}^{n}f^2_{x_i x_j}f^2_{x_i}+\sum_{\substack{i,j=1\\ i\neq j}}^{n}2f_{x_i x_i}f_{x_i x_j}f_{x_i}f_{x_j} \\& 
\leq \sum_{i,j=1}^{n}f^2_{x_i x_j}f^2_{x_i}+\sum_{\substack{i,j=1\\ i\neq j}}^{n}(f^2_{x_i x_i}f^2_{x_j}+f^2_{x_i x_j}f^2_{x_i})\\ &
 =  \sum_{i,j=1}^{n}f^2_{x_i x_i}f^2_{x_j}+2\sum_{\substack{i,j=1\\ i\neq j}}^{n}f^2_{x_i x_j}f^2_{x_i}=\lvert D^2f\rvert^2 \lvert \nabla f\rvert^2.
 \end{split}
  \]
Relations \eqref{InequalityLaplacian} and \eqref{InequalityGradienHessian} are so shown. As to \eqref{InequalityHessian}, this is Lemma 2.2 of \cite{LankeitWangConsumptLogistic}.
\end{proof}
\end{lemma}
\begin{lemma} (Gagliardo-Nirenberg inequality)\label{InequalityG-NLemma}
Let $\Omega$ be a bounded Lipschitz domain of $\R^n,$ $1\leq  \mathfrak{q}, \mathfrak{r}\leq \infty$ and  $\mathfrak{j}$ and $\mathfrak{m}$  integers satisfying $0\leq \mathfrak{j}<\mathfrak{m}$. Moreover, let be $\mathfrak{p}\in \mathbb{R}^+$ and $\frac{\mathfrak{j}}{\mathfrak{m}} \leq \theta \leq 1$  such that this equality 
$\frac{1}{\mathfrak{p}} = \frac{\mathfrak{j}}{n} + \left( \frac{1}{\mathfrak{r}} - \frac{\mathfrak{m}}{n} \right) \theta  + \frac{1 - \theta}{\mathfrak{q}}$ holds.Then there exists a constant $C_{GN}$  such that for all $f\in L^\mathfrak{q}(\Omega)$, with $D^{\mathfrak{m}} f$ in  $L^\mathfrak{r}(\Omega)$,
\begin{equation}\label{InequalityTipoG-N} 
\| D^{\mathfrak{j}} f \|_{L^{\mathfrak{p}}(\Omega)} \leq C_{GN} ( \| D^{\mathfrak{m}} f \|_{L^{\mathfrak{r}}(\Omega)}^{\theta} \| f \|_{L^{\mathfrak{q}}(\Omega)}^{1 - \theta}+  \| f \|_{L^{\mathfrak{s}}(\Omega)}),
\end{equation}
with arbitrary $\mathfrak{s}\in \mathbb{R}^+.$
\end{lemma}
We also make use of these general results. 
\begin{lemma}\label{BoundsInequalityLemmaTecnicoYoung}  
Let $A,B \geq 0$ and $d_1,d_2>0$. Then for 
\begin{equation}\label{ConstantForTechincalInequality}
k:=\min\{d_1,d_2\}\quad \textrm{and}\quad 
d_3=d_3(k):=\frac{d_1-k}{d_1}\big(\frac{d_1}{k}\big)^\frac{k}{d_1-k}+\frac{d_2-k}{d_2}\big(\frac{d_2}{k}\big)^\frac{k}{d_2-k},
\end{equation}
we have 
\begin{equation}\label{InequalityForFinallConclusion}
A^{d_1}+B^{d_2}\geq 2^{-k}(A+B)^k-d_3.
\end{equation}
\begin{proof}
Two applications of the Young inequality with exponents $\frac{d_1}{k}$ and $\frac{d_1}{d_1-k}$ and  $\frac{d_2}{k}$ and $\frac{d_2}{d_2-k}$ infer, respectively
\begin{equation*}
A^k\frac{d_1}{k} \leq \frac{k}{d_1}A^{d_1}+\frac{d_1-k}{d_1}\big(\frac{d_1}{k}\big)^\frac{d_1}{d_1-k},
\end{equation*}
and
\begin{equation*}
B^k\frac{d_2}{k} \leq \frac{k}{d_2} B^{d_2}+\frac{d_2-k}{d_2}\big(\frac{d_2}{k}\big)^\frac{d_2}{d_2-k},
\end{equation*}
which imply 
\begin{equation*}
A^k \leq \big(\frac{k}{d_1}\big)^2A^{d_1}+\frac{d_1-k}{d_1}\big(\frac{d_1}{k}\big)^\frac{k}{d_1-k}\leq  A^{d_1}+\frac{d_1-k}{d_1}\big(\frac{d_1}{k}\big)^\frac{k}{d_1-k},
\end{equation*}
and
\begin{equation*}
B^k \leq \big(\frac{k}{d_2}\big)^2B^{d_2}+\frac{d_2-k}{d_2}\big(\frac{d_2}{k}\big)^\frac{k}{d_2-k}\leq  B^{d_2}+\frac{d_2-k}{d_2}\big(\frac{d_2}{k}\big)^\frac{k}{d_2-k}.
\end{equation*}
Finally, plugging these last inequalities into 
\begin{equation}\label{AlgebraicInequality2toalpha}  
(A+B)^k\leq 2^k(A^k+B^k),
\end{equation}
valid for any $A,B\geq 0$ and $k>0$, we conclude. We have to observe that for $i=1,2$ it holds that 
\begin{equation*}
\frac{d_i-k}{d_i}\big(\frac{d_i}{k}\big)^\frac{k}{d_i-k}\rightarrow 0 \quad \textrm{as}\quad d_i\rightarrow k,
\end{equation*}
so that the expression of $d_3(k)$ in \eqref{ConstantForTechincalInequality} is meaningful for both $i=1,2$.
\end{proof}
\end{lemma}
\begin{lemma}\label{Lemmapbarra}  
For any $n\in \mathbb{N}$,  $q_1>n+2$, $q_2>\frac{n+2}{2}$ and $m,\alpha \in \R$ with $\alpha < \frac{m+1}{2}$ let
\begin{equation}\label{ConstantForTechincalInequality_Barp}
\bar{p}:=\max 
\begin{Bmatrix}
\frac{n}{2}(1-m) \vspace{0.1cm}  \\  \frac{q_1(2\alpha+1)}{2}\vspace{0.1cm}  \\ 1+m-2\alpha \vspace{0.1cm}  \\ \frac{q_1}{2} \vspace{0.1cm} \\ 1-m\frac{(n+1)q_1-(n+2)}{q_1-(n+2)} \vspace{0.1cm}  \\ 1-\frac{m}{1-\frac{n}{n+2}\frac{q_2}{q_2-1}}  
\end{Bmatrix}
+1.
\end{equation}
Then these relations hold:
\begin{equation}\label{InequalityForG-NInu^p+1} 
0< \frac{\frac{n}{2}(m+p-1)(1-\frac{1}{p})}{1-\frac{n}{2}+\frac{n}{2}(m+p-1)}<1\quad \textrm{for all} \quad p \geq \bar{p},
\end{equation}
\begin{equation}\label{InequalityForHolderInu^p+1}
0< \frac{p+2\alpha -m-1}{p}<1\quad \textrm{for all} \quad p \geq \bar{p},
\end{equation}
\begin{equation}\label{AAA} 
p>\frac{q_1}{2}\quad \textrm{for all} \quad p \geq \bar{p},
\end{equation}
\begin{equation}\label{AAA1} 
p> 1-m\frac{(n+1)q_1-(n+2)}{q_1-(n+2)}\quad \textrm{for all} \quad p \geq \bar{p}.
\end{equation}
\begin{equation}\label{AAA2} 
p>1-\frac{m}{1-\frac{n}{n+2}\frac{q_2}{q_2-1}} \quad \textrm{for all} \quad p \geq \bar{p}.
\end{equation}
\begin{equation}\label{AAA0}  
p>\frac{n}{2}(1-m)\quad \textrm{for all} \quad p \geq \bar{p},
\end{equation}
\begin{equation}\label{AAA3} 
p>\frac{q_1(2\alpha +1)}{2} \quad \textrm{for all} \quad p \geq \bar{p}.
\end{equation}
\begin{proof}
From the expression of $\bar{p}$, we have first that $\frac{m+\bar{p}-1}{2}\geq \frac{m+\bar{p}-1}{2\bar{p}}$ and, than, that $\frac{m+\bar{p}-1}{2\bar{p}}>\frac{n-2}{2n}$ since $\bar{p}>\frac{n}{2}(1-m)$; therefore $1-\frac{n}{2}+\frac{n}{2}(m+p-1)>0$ and successively also \eqref{InequalityForG-NInu^p+1} is attained. In addition, the remaining inequalities are clearly verified for any $p\geq \bar{p}$ once $\bar{p}$ is defined as in \eqref{ConstantForTechincalInequality_Barp}.
\end{proof}
\end{lemma}
\section{Existence of local-in time solutions and their properties}\label{ExistenceSolutionRegularizingSection}
Let us firstly give a result concerning local-in-time existence of classical solutions to system \eqref{problem}; its proof is obtained by well-established methods
involving standard parabolic regularity theory and an appropriate fixed point framework (see, for instance, \cite{Cieslak} and \cite{HorstWink}).
\begin{lemma}\label{LocalExistenceLemma}
Let $\Omega$ be a smooth and bounded domain of $\mathbb{R}^n$, with $n\geq 1$. For given $m,k\in \R$, and $\mu$ positive, let us assume that  $\chi\in C^2([0,\infty))$ satisfies relation \eqref{GrowthConditionOnS} for some $\alpha$ as in  \eqref{conditiononAlpha}. Then for any couple of nonnegative functions $(u_0,v_0)\in(W^{1,r}(\Omega))^2$, with $r>\max\{n,2\}$, problem \eqref{problem} admits a unique  local-in-time classical solution
\begin{equation*}
(u,v)\in (C([0,T_{max});W^{1,r}(\Omega))\cap C^{2,1}(\bar{\Omega}\times (0,T_{max})))^2,
\end{equation*}
where $T_{max}$ denotes the maximal existence time. Moreover, we have
\begin{equation}\label{MaximumPricnipleRelation} 
u\geq 0\quad 0\leq v\leq \lVert v_0\rVert_{L^\infty(\Omega)}\quad \textrm{in}\quad \Omega \times (0,T_{max}), 
\end{equation}
and if $T_{max}<\infty$ then
\begin{equation} \label{extensibility_criterion_Eq} 
\limsup_{t\nearrow T_{max}}(\lVert u (\cdot,t)\rVert_{L^\infty(\Omega)}+\lVert u (\cdot,t)\rVert_{W^{1,r}(\Omega)})=\infty.
\end{equation}
\begin{proof}
In line with the works \cite{ChoiWang} and \cite{DelgadoEtAl}, let us rewrite the initial-boundary value problem \eqref{problem} as 
\begin{equation}\label{compactProblemForLocalExistence}
\begin{cases}
w_t=\nabla \cdot (\mathcal{A}(w) \nabla w)+\mathcal{F}(w),\\ 
\frac{\partial w}{\partial \nu}=0\quad \textrm{on}\quad \partial \Omega\times [0,\infty),\\
w(\cdot,0)=(u_0,v_0)\quad \textrm{in}\quad \Omega,
\end{cases}
\end{equation}
where $w=(u,v)$ and 
\begin{equation*}
\mathcal{A}(w)=
\begin{pmatrix}
(u+1)^{m-1}&-(u+1)^\alpha \chi(v)\\
0&1
\end{pmatrix}
\quad \textrm{and}\quad
\mathcal{F}(w)=
\begin{pmatrix}
ku-\mu u^2\\
-uv
\end{pmatrix}
.
\end{equation*}
Then, Theorems 14.4 of \cite{AmannBook} warrants that problem \eqref{compactProblemForLocalExistence} possesses a maximal weak $W^{1,r}$-
solution. In turn, Theorem 14.6  of the same \cite{AmannBook} asserts that actually such a solution is classical and that the equation is indeed verified point-wise. In addition, if the extensibility criterion \eqref{extensibility_criterion_Eq} holds, Theorem 15.5 of (again) \cite{AmannBook} allows us to  conclude that $T_{max}=\infty$, namely that the solution is global. Accordingly, since $u_0\geq$ and $v_0\geq 0$ the maximum principle  and $\mathcal{F}\geq 0$ (see \cite{LSUBookInequality}) apply to yield both expressions in  \eqref{MaximumPricnipleRelation}.
\end{proof}
\end{lemma}
\begin{remark}
As it can be observed, the hypothesis for the proofs of Theorems 14.4, 14.5 and 14.6 of \cite{AmannBook} are more general that those fixed in our Lemma \ref{LocalExistenceLemma}, so that it holds also under weaker restrictions on the data of problem \eqref{problem}. Nevertheless,  for the ease of reading, we consider that in this present investigation it is more appropriate to fix also for the aforementioned initiatory Lemma  all the suitable assumptions which are necessary to the demonstration of the main Theorem \ref{MainTheorem}.
\end{remark}
\begin{lemma}
Let $\Omega$ be a smooth and bounded domain of $\mathbb{R}^n$, with $n\geq 1$. For any couple of nonnegative functions $(u_0,v_0)\in (W^{1,r}(\Omega))^2$, with $r>\max\{n,2\} $, let $(u,v)$ be the local-in-time classical solution of problem \eqref{problem} provided by Lemma \ref{LocalExistenceLemma}. Then we have
\begin{equation}\label{Bound_of_u} 
\int_\Omega u (\cdot ,t)  \le m\quad \textrm{for all}\quad  t\in (0,T_{max}),
\end{equation}
and
\begin{equation}\label{Bound_of_NablavSquare} 
\int_\Omega |\nabla v(\cdot, t)|^2\leq M\quad \textrm{for all}\quad  t\in (0,T_{max}),
\end{equation}
where 
\begin{equation*}
\begin{cases}
m =\max\{\frac{k_+|\Omega|}{\mu},\int_\Omega u_0\} \quad \textrm{with}\quad k_+=\max\{k,0\},\\
M=\max\{\lVert v_0\rVert^2_{L^\infty(\Omega)} (|\Omega|+2m+\frac{k_++1}{\mu}m),\int_\Omega \lvert \nabla v_0\rvert^2+\frac{\lVert v_0\rVert^2_{L^\infty(\Omega)}}{\mu}\int_\Omega  u_0\}.
\end{cases}
\end{equation*}
\begin{proof}
Taking into consideration the no-flux boundary conditions for problem \eqref{problem},  an integration of its first equation over $\Omega$ and an application of the Hölder inequality provide 
\begin{equation}\label{IntegrationOnOmegaFirstEq}
\frac{d}{dt}\int_\Omega u = k \int_\Omega u-\mu\int_\Omega u^2\leq k_+ \int_\Omega u-\frac{\mu}{|\Omega|}\bigg(\int_\Omega u\bigg)^2\quad \textrm{for all}\quad t\in (0,T_{max}),
\end{equation}
so that \eqref{Bound_of_u}  is a consequence of an ODI-comparison argument. 

As to \eqref{Bound_of_NablavSquare}, from the second equation of \eqref{problem},  and again through an integration over $\Omega$, the Young inequality, the bound for $v$ given in \eqref{MaximumPricnipleRelation} and \eqref{Bound_of_u} entail that for all $t\in (0,T_{max})$
\begin{equation}\label{AA}
\begin{split}
\frac{d}{dt}\int_\Omega \lvert \nabla v\rvert^2&=2\int_\Omega \nabla v \cdot \nabla (\Delta v -uv)=-2\int_\Omega (\Delta v)^2+2\int_\Omega uv \Delta v \\&
-2 \int_\Omega (\Delta v)^2+2 \int_\Omega v(u-1)\Delta v-2\int_\Omega  \lvert \nabla v\rvert^2  \\&
 \leq - \int_\Omega (\Delta v)^2-2\int_\Omega  \lvert \nabla v\rvert^2 + \int_\Omega v^2(u-1)^2\\ & \leq -\int_\Omega  \lvert \nabla v\rvert^2 + \lVert v_0\rVert^2_{L^\infty(\Omega)}\int_\Omega u^2+ 2\lVert v_0\rVert^2_{L^\infty(\Omega)}m+ \lVert v_0\rVert^2_{L^\infty(\Omega)}|\Omega|.
\end{split}
\end{equation}
On the other hand, from the equality in \eqref{IntegrationOnOmegaFirstEq}   we arrive at 
\begin{equation}\label{BB}
\begin{split}
\frac{\lVert v_0\rVert^2_{L^\infty(\Omega)}}{\mu}\frac{d}{dt}\int_\Omega u &=  \frac{k\lVert v_0\rVert^2_{L^\infty(\Omega)}}{\mu} \int_\Omega u-\lVert v_0\rVert^2_{L^\infty(\Omega)}\int_\Omega u^2\\ & \leq  \frac{k_+\lVert v_0\rVert^2_{L^\infty(\Omega)}}{\mu}\int_\Omega u-\lVert v_0\rVert^2_{L^\infty(\Omega)}\int_\Omega u^2, 
\\ &
\end{split}
\end{equation}
for all $t\in (0,T_{max})$. Thereafter, by adding \eqref{AA} to \eqref{BB} and by virtue of \eqref{Bound_of_u}, we obtain  this inequality for $h(t):=\int_\Omega \lvert \nabla v\rvert^2 + \frac{\lVert v_0\rVert^2_{L^\infty(\Omega)}}{\mu} \int_\Omega u$
\begin{equation*}
h'(t)\leq -h(t)+\bigg(2\lVert v_0\rVert^2_{L^\infty(\Omega)}+ \frac{(k_++1)\lVert v_0\rVert^2_{L^\infty(\Omega)}}{\mu}\bigg)m+\lVert v_0\rVert^2_{L^\infty(\Omega)}|\Omega|,
\end{equation*}
so that we can conclude thanks to a further comparison argument.
\end{proof}
\end{lemma}
\section{A priori estimates and proof of the main theorem}\label{EstimatesAndProofSection}
In this section our principal objective is to gain some uniform bounds for both $u$ and $v$. In particular, we aid to control with a suitable positive and time independent constant  $\lVert u\lVert_{L^p(\Omega)}$ and $\lVert \nabla v\lVert_{L^{2p}(\Omega)}$, for $p$ sufficiently large  and on the whole interval $(0,T_{max})$; this is attained by establishing an absorptive differential inequality for $\Phi(t)=\int_\Omega (u+1)^p +\chi_0^{2p}\int_\Omega |\nabla v|^{2p}$.
\begin{lemma}\label{Estim_general_For_u^pLemma} 
Let $\Omega$ be a smooth and bounded domain of $\mathbb{R}^n$, with $n\geq 1$. For any couple of nonnegative functions  $(u_0,v_0)\in (W^{1,r}(\Omega))^2$, with $r>\max\{n,2\} $, let $(u,v)$ be the local-in-time classical solution of problem \eqref{problem} provided by Lemma \ref{LocalExistenceLemma}. Then, for any $p\geq \bar{p}$,  $\bar{p}$ being the constant given in \eqref{ConstantForTechincalInequality_Barp}, and $\epsilon_1,\epsilon_2,\epsilon_3$ positive real numbers,  we have
\begin{equation}\label{Estim_general_For_u^p}
\begin{split}
&\frac{d}{dt}\int_\Omega (u+1)^p+\frac{4p(p-1)}{(m+p-1)^2}(1-\chi_0 \epsilon_1)\int_\Omega \lvert \nabla (u+1)^\frac{m+p-1}{2}\rvert^2 \leq \\ &
\quad +(\epsilon_3+\chi_0 p(p-1)C_1(\epsilon_1)C_2(\epsilon_2)-\mu p)\int_\Omega (u+1)^{p+1}\\ &
\quad+ p(p-1)\chi_0 C_1(\epsilon_1)\epsilon_2 2 (4p^2+n)\lVert v_0\Vert^2_{L^\infty(\Omega)}\int_\Omega \lvert \nabla  v\rvert^{2p-2} \lvert D^2 v\rvert^2\\ &
\quad +|\Omega|p((2\mu+k_+)C_3(\epsilon_3)+c_0(p-1)\chi_0),
\end{split}
\end{equation}
where 
\[
\begin{cases}
C_1(\epsilon_1)=\frac{1}{4 \epsilon_1}\quad C_2(\epsilon_2)=\frac{p}{p+1} (\epsilon_2(p+1))^\frac{-1}{p} \quad    C_3(\epsilon_3)=\frac{1}{p+1}(\frac{\epsilon_3(p+1)}{(2\mu+k_+) p^2})^{-p}\\
 c_0=C_1(\epsilon_1)C_2(\epsilon_2)\frac{m+1-2\alpha}{p}(\frac{p}{p+2\alpha-m-1})^\frac{p+2\alpha-m-1}{2\alpha-m-1}.
\end{cases}\]
\begin{proof}
For $\bar{p}$ as in \eqref{ConstantForTechincalInequality_Barp}, let $p\geq \bar{p}$ and $\epsilon_1,\epsilon_2,\epsilon_3$ positive real numbers which will be properly chosen through the paper. Testing the first equation of problem \eqref{problem} by $p(u+1)^{p-1}$, using its boundary conditions and relation \eqref{GrowthConditionOnS} provide
\begin{equation}\label{Estim_1_For_u^p}  
\begin{split}
\frac{d}{dt} \int_\Omega (u+1)^p&=p\int_\Omega (u+1)^{p-1}u_t \leq-p(p-1)\int_\Omega (u+1)^{p+m-3}\lvert \nabla u\rvert^2  \\ &\quad +p(p-1)\chi_0 \int_\Omega  (u+1)^{p+\alpha-2}\nabla u \cdot  \nabla v
 \\ &\quad+kp\int_\Omega u (u+1)^{p-1}-\mu p\int_\Omega u^2(u+1)^{p-1}, 
%
\end{split}
\end{equation}
on  $(0,T_{max})$. Now, the Young inequality entails for all $t\in (0,T_{max})$
\begin{equation}\label{Estim_3_For_u^p} 
\begin{split}
 \int_\Omega    (u+1)^{p+\alpha-2}\nabla u \cdot  \nabla v& \leq  \epsilon_1 \int_\Omega (u+1)^{p+m-3}|\nabla u|^2\\ &\quad 
+C_1(\epsilon_1)\int_\Omega (u+1)^{p+2\alpha-m-1}|\nabla v|^2.
\end{split}
\end{equation}
Since from \eqref{InequalityForHolderInu^p+1} we have that $0<\frac{p+2\alpha-m-1}{p}<1,$ applications of the Young inequality give on $(0,T_{max})$
\begin{equation}\label{Estim_uGradv}   
\begin{split}
C_1(\epsilon_1)\int_\Omega (u+1)^{p+2\alpha-m-1}|\nabla v|^2&\leq C_1(\epsilon_1) \epsilon_2\int_\Omega \lvert \nabla v \rvert^{2(p+1)}
\\ &\quad+C_1(\epsilon_1) C_2(\epsilon_2) \int_\Omega (u+1)^\frac{(p+1)(p+2\alpha -m-1)}{p}\\ &
\leq C_1(\epsilon_1) \epsilon_2\int_\Omega \lvert \nabla v \rvert^{2(p+1)}
\\ &\quad+C_1(\epsilon_1) C_2(\epsilon_2) \int_\Omega (u+1)^{p+1}+c_0|\Omega|.
%
\end{split}
\end{equation}
As to the contribution from the logistic source, for all $t \in (0,T_{max})$ we can write
\begin{equation}\label{EstimateLogisticPartt}   
\begin{split}
&kp\int_\Omega u (u+1)^{p-1}-\mu p\int_\Omega u^2(u+1)^{p-1}\leq\\ & k_+p\int_\Omega   (u+1)^{p}-\mu p\int_\Omega   (u+1)^{p+1} +2\mu p\int_\Omega (u+1)^{p},
\end{split}
\end{equation}
where we have employed the inequality $-u^2\leq -(u+1)^2+2(u+1).$ Successively, the Young inequality enables us to deduce that on $(0,T_{max})$
\begin{equation} 
 p(2\mu+k_+)\int_\Omega (u+1)^{p}\leq \epsilon_3\int_\Omega (u+1)^{p+1}+p(2\mu+k_+)C_3(\epsilon_3)|\Omega|.\\
\end{equation}
Taking into account  that 
\begin{equation}\label{ReducedExpressioForGradientu^m+p-1}
 p(p-1)\int_\Omega (u+1)^{p+m-3}\lvert \nabla u\rvert^2=\frac{4p(p-1)}{(m+p-1)^2}\int_\Omega |\nabla (u+1)^\frac{m+p-1}{2}|^2,
\end{equation} 
and that  \eqref{InequalityHessian}, with $q=p$, and the bound for $v$ in $\Omega \times (0,T_{max})$ in \eqref{MaximumPricnipleRelation} imply
\begin{equation}\label{ExpressioForGradientv^2p+2}
\int_\Omega \lvert \nabla u\rvert^{2p+2}\leq2 (4p^2+n)\lVert v_0\Vert^2_{L^\infty(\Omega)}\int_\Omega \lvert \nabla  v\rvert^{2p-2} \lvert D^2 v\rvert^2,
\end{equation} 
our thesis is justified once \eqref{Estim_1_For_u^p}-\eqref{ExpressioForGradientv^2p+2} are collected. 
\end{proof}
\end{lemma}
The main idea of the following lemma is not new but comes from Lemma 4.2 of \cite{LankeitWangConsumptLogistic}; despite, we have to adapt those steps to our scope. 
\begin{lemma}\label{Estim_general_For_nablav^2pLemma} 
Let $\Omega$ be a smooth and bounded domain of $\mathbb{R}^n$, with $n\geq 1$. For any couple of nonnegative functions $(u_0,v_0)\in (W^{1,r}(\Omega))^2$, with $r>\max\{n,2\} $, let $(u,v)$ be the local-in-time classical solution of problem \eqref{problem} provided by Lemma \ref{LocalExistenceLemma}. Then, for any $p\geq 1$ and $\delta_1>0$,  we have that  for all $t \in (0,T_{max})$
\begin{equation}\label{Estim_general_For_Nablav^2p}
\begin{split}
&\frac{d}{dt}\int_\Omega  |\nabla v|^{2p}+2p(1- (p+n-1)(4p^2+n)\delta_1\lVert v_0 \rVert^4_{L^\infty(\Omega)}) \int_\Omega |\nabla v|^{2p-2} |D^2v|^2 \leq \\ &
\quad +p(p+n-1)D_1(\delta_1)\lVert v_0 \rVert^2_{L^\infty(\Omega)}\int_\Omega (u+1)^{p+1},
\end{split}
\end{equation}
where $D_1(\delta_1)=\frac{2}{p+1}(\delta_1\frac{p+1}{p-1})^\frac{1-p}{2}$.
\begin{proof}
From the second equation of \eqref{problem}, we derive this equality valid for all $x\in \Omega$ and $t\in(0,T_{max})$:
\begin{equation*}
\begin{split}
(|\nabla v|^2)_t&=2 \nabla v\cdot \nabla v_t=2 \nabla v \cdot \nabla \Delta v -2\nabla v \cdot \nabla(uv)\\ &
= \Delta |\nabla v|^2 -2 |D^2 v|^2  -2\nabla v \cdot \nabla (uv).
\end{split}
\end{equation*}
Successively, multiplying this last relation by $\lvert \nabla v\rvert^{2p-2}$ and integrating over $\Omega$ lead to
\begin{equation}\label{A_00}
\begin{split}
&\frac{1}{p}\frac{d}{dt}\int_\Omega  |\nabla v|^{2p}+(p-1)\int_\Omega |\nabla v|^{2p-4}|\nabla |\nabla v|^2|^2+2 \int_\Omega |\nabla v|^{2p-2} |D^2v|^2\\ &
\quad \leq -2 \int_\Omega |\nabla v|^{2p-2} \nabla v \cdot \nabla (uv)\quad \textrm{for all}\quad t \in (0,T_{max}).
\end{split}
\end{equation}
Hence, an integration by parts to the right hand side term produces, also thanks again to \eqref{MaximumPricnipleRelation},
\begin{equation}\label{A_01} 
\begin{split}
& -2 \int_\Omega |\nabla v|^{2p-2} \nabla v \cdot \nabla (uv)= 2 \int_\Omega uv|\nabla v|^{2p-2}\Delta v  \\& \quad   +2(p-1) \int_\Omega uv|\nabla v|^{2p-4}\nabla v \cdot \nabla \lvert\nabla v\rvert^2 \\& 
 \leq 2 \lVert v_0\rVert_{L^\infty(\Omega)}\int_\Omega u|\nabla v|^{2p-2}|\Delta v|\\&\quad +2 (p-1) \lVert v_0\rVert_{L^\infty(\Omega)}\int_\Omega u|\nabla v|^{2p-3}\lvert\nabla \lvert \nabla v\rvert^2\rvert \quad \textrm{for all}\quad t \in (0,T_{max}).
\end{split}
\end{equation}
In addition, the Young and \eqref{InequalityLaplacian} inequalities allow us to derive
\begin{equation}\label{Estimateu^2gradv^2_1}
\begin{split}
& 2 \lVert v_0\rVert_{L^\infty(\Omega)}\int_\Omega u|\nabla v|^{2p-2}|\Delta v|\leq \frac{1}{n}\int_\Omega |\nabla v|^{2p-2}|\Delta v|^2\\ & 
\quad +n\lVert v_0\rVert_{L^\infty(\Omega)}^2\int_\Omega u^2|\nabla v|^{2p-2}\leq  \int_\Omega |\nabla v|^{2p-2}|D^2 v|^2\\ & 
\quad +n\lVert v_0\rVert_{L^\infty(\Omega)}^2\int_\Omega u^2|\nabla v|^{2p-2}\quad \textrm{for all}\quad t \in (0,T_{max}),
 \end{split}
\end{equation}
and similarly
\begin{equation}\label{Estimateu^2gradv^2_Second}
\begin{split}
& 2(p-1) \lVert v_0\rVert_{L^\infty(\Omega)}\int_\Omega u |\nabla v|^{2p-3}\lvert\nabla \lvert \nabla v\rvert^2\rvert\leq (p-1)\int_\Omega |\nabla v|^{2p-4}\lvert \nabla \lvert \nabla v\rvert^2\rvert^2\\ & 
\quad +\lVert v_0\rVert_{L^\infty(\Omega)}^2(p-1)\int_\Omega u^2|\nabla v|^{2p-2}\quad \textrm{for all}\quad t \in (0,T_{max}).
 \end{split}
\end{equation}
Now, for $\delta_1>0$, the Young inequality with exponents $\frac{p+1}{p-1}$ and $\frac{p+1}{2}$ produces on $(0,T_{max})$ also
\begin{equation}\label{Estimatingnablav^2p+2}
\begin{split}
\int_\Omega u^2|\nabla v|^{2p-2}&\leq \delta_1 \int_\Omega |\nabla v|^{2p+2}+D_1(\delta_1)\int_\Omega u^{p+1} \\ & \leq \delta_1 \int_\Omega |\nabla v|^{2p+2}+D_1(\delta_1)\int_\Omega (u+1)^{p+1}.
 \end{split}
\end{equation}
By virtue of \eqref{ExpressioForGradientv^2p+2}, we have the claim introducing \eqref{A_01}-\eqref{Estimatingnablav^2p+2} into \eqref{A_00}.
\end{proof}
\end{lemma}
\begin{lemma}\label{LemmaAbsorptiveMainInequality}
Let $\Omega$ be a smooth and bounded domain of $\mathbb{R}^n$, with $n\geq 1$. For any couple of nonnegative functions $(u_0,v_0)\in (W^{1,r}(\Omega))^2$, with $r>\max\{n,2\} $, let $(u,v)$ be the local-in-time classical solution of problem \eqref{problem} provided by Lemma \ref{LocalExistenceLemma}. Let also assume that for any $p\geq \bar{p}$,  where $\bar{p}$ is the constant given by \eqref{ConstantForTechincalInequality_Barp},  
$\mu$ satisfies the following relation
\begin{equation}\label{AssumptionOnMu}
\mu \geq k_1(n,p)\lVert \chi_0 v_0 \rVert_{L^\infty(\Omega)}^\frac{2}{p}+k_2(n,p)\lVert \chi_0 v_0 \rVert_{L^\infty(\Omega)}^{2p},
\end{equation}
where 
\begin{equation*}
\begin{cases}
k_1(p,n)=p^2\big(\frac{p-1}{p+1}\big)^\frac{p+1}{p}(4p^2+n)^\frac{1}{p}(\chi_0\lVert v_0\rVert_{L^\infty(\Omega)})^\frac{2}{p},\\
k_2(p,n)=\frac{p}{p+1}2^p(p+n-1)^\frac{p+1}{2}\big(\frac{p-1}{p+1}\big)^\frac{p-1}{2}(4p^2+n)^\frac{p-1}{2}.
\end{cases}
\end{equation*}
Then  there exists a positive constant $L_1$ such that for any $p\geq \bar{p}$
\begin{equation}\label{BoundU+1EnadNabla}
\int_\Omega (u+1)^p +\chi_0^{2p}\int_\Omega |\nabla v|^{2p}\leq L_1 \quad \textrm{for all}\quad t \in (0,T_{max}).
\end{equation}
\begin{proof}
For this particular choice of the constants $\epsilon_i$ ($i=1,2,3$) and $\delta_1$ introduced above,
\begin{equation*}
\begin{cases}
\epsilon_1=\frac{1}{2\chi_0}, \quad \epsilon_2=\frac{\chi_0^{2p-1}}{4(p-1)C_1(\epsilon_1)(4p^2+n)\lVert v_0\rVert_{L^\infty(\Omega)}^2},\quad \\
\delta_1=\frac{1}{4(p+n-1)(4p^2+n)\lVert v_0\rVert_{L^\infty(\Omega)}^4},\\
 \epsilon_3=\frac{p^2}{2}\big(\frac{p-1}{p+1}\big)^\frac{p+1}{p}(4p^2+n)^\frac{1}{p}(\chi_0\lVert v_0\rVert_{L^\infty(\Omega)})^\frac{2}{p},
\end{cases}
\end{equation*}
let us multiply expression \eqref{Estim_general_For_Nablav^2p} by $\chi_0^{2p}$ and, then, let us add the result to relation \eqref{Estim_general_For_u^p}. Hence, for $\Phi(t)=\int_\Omega (u+1)^p +\chi_0^{2p}\int_\Omega |\nabla v|^{2p}$, we have
\begin{equation*}\label{Estim_general_For_Phi_1}
\begin{split}
&\Phi'(t)+\frac{2p(p-1)}{(m+p-1)^2}\int_\Omega \lvert \nabla (u+1)^\frac{m+p-1}{2}\rvert^2+p\chi_0^{2p} \int_\Omega |\nabla v|^{2p-2} |D^2v|^2 \leq \\ &
 +p( k_1(n,p)\lVert \chi_0 v_0 \rVert_{L^\infty(\Omega)}^\frac{2}{p}+k_2(n,p)\lVert \chi_0 v_0 \rVert_{L^\infty(\Omega)}^{2p}-\mu)\int_\Omega (u+1)^{p+1}\\ &
+|\Omega|p((2\mu+k_+)C_3(\epsilon_3)+c_0(p-1)\chi_0),
\end{split}
\end{equation*}
which in view of assumption \eqref{AssumptionOnMu} reads
\begin{equation}\label{Estim_general_For_Phi_2}
\begin{split}
&\Phi'(t)+\frac{2p(p-1)}{(m+p-1)^2}\int_\Omega \lvert \nabla (u+1)^\frac{m+p-1}{2}\rvert^2+\frac{\chi_0^{2p}}{p} \int_\Omega \vert \nabla \lvert \nabla v\rvert^p\rvert^2 \leq c_1,
\end{split}
\end{equation}
where $c_1=|\Omega|p((2\mu+k_+)C_3(\epsilon_3)+c_0(p-1)\chi_0)$ and where, thanks to \eqref{InequalityGradienHessian} we have employed
\[\vert \nabla \lvert \nabla v\rvert^p\rvert^2=\frac{p^2}{4}\lvert \nabla v \rvert^{2p-4}\vert \nabla \lvert \nabla v\rvert^2\rvert^2=p^2\lvert \nabla v \rvert^{2p-4}\lvert D^2v \nabla v \rvert^2\leq p^2|\nabla v|^{2p-2} |D^2v|^2.\] 
Now, the Gagliardo-Nirenberg inequality \eqref{InequalityTipoG-N} with $\mathfrak{j}=0$, $\mathfrak{p}=\frac{2p}{m+p-1}$, $\mathfrak{m}=1$, $\mathfrak{r}=2$ and $\mathfrak{s}=\mathfrak{q}=\frac{2}{m+p-1}$ infers in conjunction to \eqref{InequalityForG-NInu^p+1}  of Lemma \ref{Lemmapbarra} that 
\[0<\theta_1=\frac{n\frac{m+p-1}{2}(1-\frac{1}{p})}{1-\frac{n}{2}+n\frac{m+p-1}{2}}<1.\]
Subsequently, we get
\begin{equation*} 
\begin{split}
\int_\Omega (u+1)^{p}&=\lvert \lvert (u+1)^\frac{m+p-1}{2}\lvert \lvert_{L^\frac{2p}{m+p-1}(\Omega)}^\frac{2p}{m+p-1} 
\\&
\leq c_2        \lvert \lvert\nabla (u+1)^\frac{m+p-1}{2}\lvert \lvert_{L^2(\Omega)}^{\frac{2p}{m+p-1}\theta_1}      \lvert \lvert (u+1)^\frac{m+p-1}{2}\lvert \lvert_{L^\frac{2}{m+p-1}(\Omega)}^{(1-\theta_1)\frac{2p}{m+p-1}} \\ &
\quad +c_2      \lvert \lvert (u+1)^\frac{m+p-1}{2}\lvert \lvert^\frac{2p}{m+p-1}_{L^\frac{2}{m+p-1}(\Omega)},
 \end{split}
\end{equation*}
where $c_2= (2C_{GN})^\frac{2p}{m+p-1}$ and having also taken into consideration  \eqref{AlgebraicInequality2toalpha}. Hence, recalling bound \eqref{Bound_of_u} and introducing $c_3=c_2 \max\{(m+|\Omega|)^{(1-\theta_1)p},(m+|\Omega|)^p\}$ the last inequality entails
\begin{equation}\label{Estim_u+1topGaglNiren}
\begin{split}
\int_\Omega (u+1)^{p}\leq& c_3\Big(\int_\Omega \lvert \nabla (u+1)^\frac{m+p-1}{2}\rvert^2\Big)^\frac{p\theta_1}{m+p-1}+c_3 \quad t \in (0,T_{max}).
 \end{split}
\end{equation}
In a similar way, again the Gagliardo-Nirenberg inequality \eqref{InequalityTipoG-N} allows us to write, for  $\mathfrak{j}=0$, $\mathfrak{p}=2$, $\mathfrak{m}=1$, $\mathfrak{r}=2$ and $\mathfrak{s}=\mathfrak{q}=\frac{2}{p}$ 
\begin{equation*}
\begin{split}
\int_\Omega \lvert \nabla v\rvert^{2p}&=\lvert \lvert \lvert \nabla v\rvert^p\lvert \lvert_{L^2(\Omega)}^2 
\\&
\leq c_4        \lvert \lvert\nabla  \lvert \nabla v \rvert^p\rvert \lvert_{L^2(\Omega)}^{2\theta_2}      \lvert \lvert\lvert \nabla v \rvert^p\lvert \lvert_{L^\frac{2}{p}(\Omega)}^{2(1-\theta_2)} 
+c_4      \lvert \lvert \lvert \nabla v \rvert^p\lvert \lvert^2_{L^\frac{2}{p}(\Omega)},
 \end{split}
\end{equation*}
where $c_4= (2C_{GN})^2$ and $0<\theta_2= \frac{\frac{np}{2}-\frac{n}{2}}{1+\frac{np}{2}-\frac{n}{2}}<1$. Successively, we have 
\begin{equation}\label{Estim_Nabla nabla v^p^2}
\begin{split}
\int_\Omega \lvert \nabla v\rvert^{2p}\leq c_5 \Big(\int_\Omega \lvert \nabla \lvert \nabla v \rvert^p\rvert^2\Big)^{\theta_2}+c_5 \quad t \in (0,T_{max}),
 \end{split}
\end{equation}
with $c_5=c_4\max\{M^{(1-\theta_2)p},M^p\}$, $M$ being the constant provided by \eqref{Bound_of_NablavSquare}.

As a consequence of all of the above, by making first use of inequality \eqref{AlgebraicInequality2toalpha}  in \eqref{Estim_u+1topGaglNiren} and \eqref{Estim_Nabla nabla v^p^2} and then inserting both results into \eqref{Estim_general_For_Phi_2}, we obtain thanks also to \eqref{InequalityForFinallConclusion} that $\Phi$ verifies this initial problem
\begin{equation*}\label{MainInitialProblemWithM}
\begin{cases}
\Phi'(t)\leq c_6-c_7 \Phi^\delta(t)\quad t \in (0,T_{max}),\\
\Phi(0)=\int_\Omega (u_0+1)^p+\chi_0^{2p}\int_\Omega |\nabla v_0|^{2p}, 
\end{cases}
\end{equation*}
with 
\[
\begin{cases}
\delta=\min\{\frac{m+p-1}{p \theta_1},\frac{1}{\theta_2}\}, \quad c_6=c_1+\frac{\chi_0^{2p}}{p}+\frac{2p(p-1)}{(m+p-1)^2}+d_3(\delta),\\ 
 c_7=2^{-\delta}\min\{\frac{2p(p-1)}{(m+p-1)^2}(2c_3)^{-\frac{m+p-1}{p\theta_1}},\frac{1}{p}(2c_5)^{-\frac{1}{\theta_2}}\}. 
\end{cases}  
\]
Consequently,  again an application of an ODE  comparison principle implies that $\Phi(t)\leq \max\{\Phi(0),\big(\frac{c_6}{c_7}\big)^\frac{1}{\delta}\}:=L_1$ for all $t\in(0,T_{max})$,  
and we conclude. 
\end{proof}
\end{lemma}
After these preparations, the proof of our main result consists in organizing the above statements and other facts as follows.
\subsubsection*{Proof of Theorem \ref{MainTheorem}}  
Let $\Omega$ be a smooth and bounded domain of $\mathbb{R}^n$, with $n\geq 1$. For given $m,k\in \R$, $\mu$ positive and $\chi\in C^2([0,\infty))$ satisfying relation \eqref{GrowthConditionOnS} for some $\alpha$ as in  \eqref{conditiononAlpha}, let  $(u,v)$ be the local-in-time classical solution of problem \eqref{problem} emanating from any couple of nonnegative functions $(u_0,v_0)\in(W^{1,r}(\Omega))^2$, whose existence is ensured by Lemma \ref{LocalExistenceLemma}. For $\bar{p}$ defined  in \eqref{ConstantForTechincalInequality_Barp}  of Lemma \ref{Lemmapbarra}, let us set $K_1(n,m,\alpha)=k_1(\bar{p},n)$ and $K_2(n,m,\alpha)=k_2(\bar{p},n)$, where $k_1(p,n)$ and $k_2(p,n)$ have been introduced in Lemma 
\ref{LemmaAbsorptiveMainInequality}; since \eqref{LargenessAssumptionMu} is satisfied, we have by continuity reasons that there exists $p>\bar{p}$ such that 
\begin{equation*}
\mu> k_1(p,n)\lVert \chi_0 v_0 \rVert_{L^\infty(\Omega)}^\frac{2}{p}+k_2(p,n)\lVert \chi_0 v_0 \rVert_{L^\infty(\Omega)}^{2p}.
\end{equation*}
Subsequently, assumption \eqref{AssumptionOnMu} holds so that relation \eqref{BoundU+1EnadNabla}  implies that 
\[u\in L^\infty((0,T_{max});L^p(\Omega)).\] Hereafter, coherently to the nomenclature  used by Tao and Winkler, the solution $u$ of system \eqref{problem} provided by Lemma \ref{LocalExistenceLemma} also classically solves  in $\Omega \times (0,T_{max})$ problem (A.1)  of Appendix A of \cite{TaoWinkParaPara} with 
\begin{equation*}
D(x,t,u)=(u+1)^{m-1},\quad f(x,t)=-(u+1)^\alpha \chi(v)\nabla v,\quad g(x,t)=\frac{k^2}{4 \mu}. 
\end{equation*}
Hence we deduce that (A.2)-(A.5), the second of (A-6) for any choice of $q_2$ and (A-7) with $p_0=p$ are as well verified on $(0,T_{max})$.  As to the first condition of (A-6), relation \eqref{AAA} allows us to apply  the Hölder inequality with exponents $\frac{q_1}{2p}$  and $\frac{2p-q_1}{2p}$; this, in conjunction with \eqref{GrowthConditionOnS},  \eqref{BoundU+1EnadNabla} and the fact that $\frac{2\alpha  p q_1}{2p-q_1}<p$ from \eqref{AAA3}, show that this bound holds on $(0,T_{max})$
\begin{equation*}
\begin{split}
\int_\Omega |f|^{q_1}&=\int_\Omega (u+1)^{\alpha q_1}\lvert  \chi (v)\rvert^{q_1}\lvert \nabla v\rvert^{q_1}\\ & \leq \chi_0^{q_1} \bigg(\int_\Omega \lvert \nabla v\rvert^{2p} \bigg)^\frac{q_1}{2p} \bigg(\int_\Omega (u+1)^\frac{2\alpha  p q_1}{2p-q_1}\bigg)^\frac{2p-q_1}{2p},
\end{split}
\end{equation*}
so that $f\in L^\infty((0,T_{max}); L^{q_1}(\Omega))$, and $q_1>n+2$. Moreover, by virtue of \eqref{AAA1}, \eqref{AAA2} and \eqref{AAA0}, also (A.8), (A.9) for $q_2>\frac{n+2}{2}$, and (A-10) of Lemma A.1. of \cite{TaoWinkParaPara} are valid, so we get for some $L_2>0$
\begin{equation}\label{BoundU+1Top}
\lVert u(\cdot,t)\rVert _{L^\infty(\Omega)}\leq L_2 \quad \textrm{for all}\quad t \in (0,T_{max}).
\end{equation}
Concerning the $v$-component, by means of the representation formula, we have
\begin{equation}\label{RepresentationForv}
v(\cdot,t) =e^{t\Delta}v_0 -\int_0^t e^{(t-s)\Delta}u(\cdot,s)v(\cdot,s)ds \quad \textrm{for all} \quad t\in (0,T_{max}).
\end{equation}
Now we invoke standard estimates for the Neumann heat semigroup (see Lemma 1.3 of \cite{WinklAggre}) which warrant the existence of positive constants $C_S$ and $\mu_1$ such that for all $t>0$ and $p\leq q\leq \infty$ 
\begin{equation}\label{LpLqEstimateGradient0}   
\lVert \nabla  e^{t\Delta} f  \lVert_{L^q(\Omega)}\leq C_S (1+t^{-\frac{1}{2}-\frac{n}{2}(\frac{1}{p}-\frac{1}{q})})e^{-\mu_1 t} \lVert f\lVert_{L^p(\Omega)}\quad \textrm{for all} \quad f\in L^p(\Omega),
\end{equation}
and for all $t>0$ and $2\leq p< \infty$ 
\begin{equation}\label{LpLqEstimateW1p}   
\lVert \nabla  e^{t\Delta} f  \lVert_{L^p(\Omega)}\leq C_S e^{-\mu_1 t} \lVert \nabla f\lVert_{L^p(\Omega)}\quad \textrm{for all} \quad f\in W^{1,p}(\Omega).
\end{equation}
Thereafter, from \eqref{RepresentationForv}, relying on \eqref{BoundU+1Top}  and the second bound in \eqref{MaximumPricnipleRelation}, we have for $r>\max\{n,2\}$ and on $(0,T_{max})$ 
\begin{equation*}
\begin{split}
&\lVert \nabla v(\cdot,t)\rVert_{L^r(\Omega)} \leq \lVert  \nabla e^{t\Delta}v_0 \rVert_{L^r(\Omega)} +\int_0^t \rVert \nabla e^{(t-s)\Delta}u(\cdot,s)v(\cdot,s) \rVert_{L^r(\Omega)}ds\\ &
 \leq C_S \lVert  \nabla v_0 \rVert_{L^r(\Omega)} +C_SL_2\lVert  v_0 \rVert_{L^\infty(\Omega)}\int_0^t (1+(t-s)^{-\frac{1}{2}})e^{-\mu_1(t-s)}ds,
\end{split}
\end{equation*}
where we have applied \eqref{LpLqEstimateGradient0} for $p=q=\infty$ and \eqref{LpLqEstimateW1p} for $p=r$. Subsequently, with the introduction of the Gamma function  $\Gamma$ we get for $t\in (0,T_{max})$
\begin{equation*}\label{Bound_v_1-q}
\begin{split}
\lVert  \nabla v (\cdot, t)\lVert_{L^r(\Omega)} \leq  c_8&=C_S\bigg[\lVert  \nabla v_0 \rVert_{L^r(\Omega)} +L_2\lVert  v_0 \rVert_{L^\infty(\Omega)}\bigg(\mu_1^{-1}+\mu_1^{-\frac{1}{2}}\Gamma \big(\frac{1}{2}\big)\bigg)\bigg].
\end{split}
\end{equation*}
This last inequality, in conjunction with \eqref{BoundU+1Top} and  the uniform bound for $v$ in \eqref{MaximumPricnipleRelation}, yield the boundedness for $\lVert u(\cdot,t)\rVert_{L^\infty(\Omega)}+\lVert v(\cdot,t)\rVert_{W^{1,r}(\Omega)}$ on $(0,T_{max})$. In turn, 
the extensibility criterion \eqref{extensibility_criterion_Eq} of Lemma \ref{LocalExistenceLemma} shows that $T_{max} = \infty$. Finally, the independence of the obtained estimates with respect to $t\in (0, T_{max}) = (0,\infty)$ establishes \eqref{GlobalBoundednessuinftyvW1infty} for a proper choice of $C$.
\qed
\subsubsection*{Acknowledgments}
The authors are members of the Gruppo Nazionale per l'Analisi Matematica, la Probabilit\`a e le loro Applicazioni (GNAMPA) of the Istituto Na\-zio\-na\-le di Alta Matematica (INdAM). G. Viglialoro gratefully acknowledges the Italian Ministry of Education, University and Research (MIUR) for the financial support of Scientific Project ``Smart Cities and Communities and Social Innovation - ILEARNTV  anywhere, anytime - SCN$\_$00307''.
\bibliography{Bibliography}{}
\bibliographystyle{abbrv}
\end{document}